\documentclass[12pt,twoside]{article}
\usepackage{graphicx}
\usepackage{hyperref}
\usepackage{amssymb}
\usepackage{amsfonts}

\date{}

\begin{document}


\centerline {\Large{\bf On systems of nonlinear  equations }}

\centerline{}



\centerline{\bf {$^1$M. Eshaghi Gordji, $^2$A. Ebadian,  $^3$M. B.
Ghaemi and $^4$J.
 Shokri}}

\centerline{$^1$ Department of Mathematics, Semnan University,
Semnan, Iran}

\centerline{$^{2,4}$Department of Mathematics, Urmia University,
Urmia, Iran}

\centerline{ $ ^3$ Department of Mathematics, Iran  University of
Science and Technology, Tehran,  Iran}

\centerline{E-mail: madjid.eshaghi@gmail.com, a.ebadian@urmia.ac.ir,
mghaemi@iust.ac.ir,  } \centerline{j.shokri@urmia.ac.ir}

\centerline{}







\newtheorem{Theorem}{\quad Theorem}[section]

\newtheorem{Definition}[Theorem]{\quad Definition}

\newtheorem{Corollary}[Theorem]{\quad Corollary}

\newtheorem{theo}{\quad Theorem}[section]
\newtheorem{lemm}[Theorem]{Lemma}
\newcommand{\un}{\underline}
\newcommand{\ov}{\overline}
\newtheorem{exam}[Theorem]{Example}

\begin{abstract}
In this paper, we introduce an iterative numerical method to solve
systems of nonlinear equations. The third-order convergence of this
method is analyzed. Several examples are given to illustrate the
efficiency of the proposed method.

\end{abstract}



\section{Introduction}\label{sec:1}

~~Let us consider the problem of finding a real zero of the
nonlinear system $F(x)=0$ which  $F:\Re^n\to \Re^n$. As notation
throughout this paper,  $\alpha\in \Re^n$ will denote the true
solution of the nonlinear system $F(x)=0$. More precisely Newton's
method may has used as the approximation of the following indefinite
integral, arising from Newton's theorem \cite{dennis},
\begin{equation}\label{eq:1}
f(x)=f(x_n)+\int_{x_n}^x f'(t)dt,
\end{equation}
 for nonlinear equation $f(x)=0$. Noor\cite{Noor} by using the combination of
midpoint quadrature rule and Trapezoidal rule for integral
(\ref{eq:1}) has introduced following iterative process for solving
$f(x)=0$,
\begin{equation}\label{eq:2}
    x_{n+1}=x_n-\frac{4f(x_n)}{f'(x_n)+2f'(\frac{x_n+y_{n}}{2})+f'(y_n)},
\end{equation}
where
\[y_n=x_n-\frac{f(x_n)}{f'(x_n)}.
\]

 Now, corresponding to (\ref{eq:1}), for nonlinear system $F(x)=0$ is
written,  Ortega\cite{ort1,ort2}:

\[
F(x)=F(x_n)+\int_{x_n}^x F'(t)dt,
\]
then we can extend the discussion to  solve system of nonlinear
equations $F(x)=0$, so  similar to (\ref{eq:2}), the following
iterative process for solving $F(x)=0$ is obtain as,
\[
x_{n+1}=x_n-4\Big[F'(x_n)+2F'\Big(\frac{x_n+y_n}{2}\Big)+F'(y_n)\Big]^{-1}F(x_n),\quad
n=0,1,\dots,
\]
 where
\[
y_n=x_n-F'(x_n)^{-1}F(x_n).
\]
where $F'(x_n)^{-1}$ is the Jaccobian Matrix of the function $F$
evaluated in $x_k$.
 we call this iterative process \emph{Midpoint-Trapezoidal Newton's
method} (MTN). In this paper, we analyze (MTN) in details  and prove
its third-order convergent theorem. Also, we have comparisons with
some other variants of Newton's method by numerical examples.

\section{Description of the methods}\label{sec:2}

 Let $F:\Omega\subseteq\Re^n\to\Re^n$ be sufficiently differentiable
 function and $\alpha$ be a zero of the system of nonlinear equations
 $F(x)=0$. From (\ref{eq:2}) as following
\begin{equation}\label{eq:3}
F(x)=F(x_n)+\int_{x_n}^xF'(t)dt,
\end{equation}
we saw in the previous section that by using rectangular rule for
above integral, classical Newton's method (CN) is obtained as
following
\[
x_{n+1}=x_n-F'(x_n)^{-1} F(x_n).
\]
where $F'(x)^{-1}$ is the Jacobian matrix of the function $F$
evaluated in $x_n$. If an estimation of (\ref{eq:3}) is made by
means of the trapezoidal rule and $x=\alpha$ is taken, then
\[
0\approx F(x_n)+\frac{1}{2}[F(x_n)+F(\alpha)](\alpha-x_n),
\]
is obtained and a new approximation $x_{n+1}$ of $\alpha$ is given
by
\[
x_{n+1}=x_n-2[F'(x_n)+F'(x_{n+1})]^{-1}F(x_n).
\]
For solving of the implicit form problem that this equation involve,
we use the $(n+1)$th  approximation of Newton method in right side,
\[
x_{n+1}=x_n-2 [F'(x_n)+F'(y_n)]^{-1}F(x_n),\quad n=0,1,\dots,
\]
where
\[
y_n=x_n-F'(x_n)^{-1}F(x_n).
\]
This iterative method will be called \emph{Trapezoidal Newton's
method} (TN). By using harmonic mean in (\ref{eq:3})and $x=\alpha$
is taken,we have
\[
0\approx
F(x_n)+\frac{2F(\alpha)F(x_n)}{F(x_n)+F(\alpha)}(\alpha-x_n),
\]
and the following iterative approximation is obtained
\[
x_{n+1}=x_n-\frac{1}{2}F(y_n)^{-1}F(x_n)^{-1}[F(x_n)+F(y_n)]F(x_n),\quad
n=0,1,\dots,
\]
where is,
 \[
 y_n=x_n-F'(x_n)^{-1}F(x_n),
\]
this variant of Newton's method is called \emph{Harmonic Newton's
method} (HN). If the midpoint rule is used to estimate integral
(\ref{eq:3}) and $x=\alpha$ is taken, it is obtained one
\[
0\approx F(x_n)+F\Big(\frac{x_n+\alpha}{2}\Big)(\alpha-x_n),
\]
then, by a approximation $x_{n+1}$ f $\alpha$,
\[
x_{n+1}=x_n-F'\Big(\frac{x_n+x_{n+1}}{2}\Big)^{-1}F(x_n),
\]

so, an alternative of Newton's method is obtained as following
\[
x_{n+1}=x_n-F'\Big(\frac{x_n+y_n}{2}\Big)^{-1}F(x_n),\quad
n=0,1,\dots,
\]
where
\[
y_n=x_n-F'(x_n)^{-1}F(x_n),
\]
this variant of Newton's method is called \emph{Midpoint Newton's
method} (MN).

Now, if the integral (\ref{eq:3}) is estimated using the combination
of midpoint quadrature rule and Trapezoidal rule and by considering
$x=\alpha$, we have
\[
0\approx
F(x_n)+\frac{1}{4}\Big[F'(x_n)+2F'\Big(\frac{x_n+\alpha}{2}\Big)+F'(\alpha)\Big](\alpha-x_n),
\]
so, a new approximation $x_{n+1}$ of $\alpha$ is concluded as
following:
\[
x_{n+1}=x_n-4\Big[F'(x_n)+2F'\Big(\frac{x_n+\alpha}{2}\Big)+F'(\alpha)\Big]^{-1}F(x_n),
\]
by using again the $(n+1)$th iteration of Newton's method in the
right side of this equation, the implicit problem is avoided. Then
\begin{equation}\label{eq:4}
x_{n+1}=x_n-4\Big[F'(x_n)+2F'\Big(\frac{x_n+y_n}{2}\Big)+F'(y_n)\Big]^{-1}F(x_n),\quad
n=0,1,\dots,
\end{equation}
is deduced, where
\[
y_n=x_n-F'(x_n)^{-1}F(x_n).
\]
This iterative process is called \emph{Midpoint-Trapezoidal Newton's
method} (MTN).

In the next section we prove that, (MTN) has third-order
convergence. The convergence of  the other variants of Newton's
methods can be proved analogously.
\section{Main result}
~~ In this section the third-order convergence of
Midpoint-Trapezoidal Newton's method (MTN) is proven by following
theorem.
\begin{theo}\label{theo:1}
Let $F:\Omega\subseteq \Re^n\to \Re^n$, is $k$-times Fr\'echet
differentiable in convex $\Omega$ containing  the root $\alpha$ of
$F(x)=0$. The Midpoint-Trapezoidal Newton's method (\ref{eq:4}) has
third-order convergence.
\end{theo}

{\bf Proof:} The Taylor's expansion for any $x,x_n\in \Omega$,
\cite{ort2}:
\[
\begin{array}{ll}
F(x)=F(x_n)+F'(x_n)(x-x_n)+\frac{1}{2!}F''(x_n)(x-x_n)^2\\
\hspace{1.3cm}+\frac{1}{3!}F^{(3)}(x-x_n)^3\newline+\cdots+\frac{1}{k!}F^{(k)}(x_n)(x-x_n)^k+\cdots,
\end{array}
\]
with $x=\alpha$ and defining $e_n=x_n-\alpha$ we have:
\[
\begin{array}{ll}
F(\alpha)=F(x_n)+F'(x_n)(\alpha-x_n)+\frac{1}{2!}F''(x_n)(\alpha-x_n)^2\\
\hspace{1.3cm}+\frac{1}{3!}F^{(3)}(\alpha-x_n)^3+\cdots+\frac{1}{k!}F^{(k)}(x_n)(\alpha-x_n)^k+\cdots,\\
\hspace{1.08cm}=F(x_n)-F'(x_n)e_n+\frac{1}{2!}F''(x_n)(\alpha-x_n)^2\\
\hspace{1.3cm}-\frac{1}{3!}F^{(3)}e_n^3+\cdots+(-1)^k\frac{1}{k!}F^{(k)}(x_n)e_n^k+\cdots
\end{array}
\]
For $k=3$ and from $F(\alpha)=0$ we have:
\begin{equation}\label{eq:5}
F(x_n)=F'(x_n)e_n-\frac{1}{2!}F''(x_n)e_n^2
+\frac{1}{3!}F^{(3)}(x_n)e_n^3+O(\|e_n\|^4).
\end{equation}
From (\ref{eq:5}) we can write the product $F'(x_n)^{-1}F(x_n)$ as
following:
\[
F'(x_n)^{-1}F(x_n)=F'(x_n)^{-1}\Big(F'(x_n)e_n-\frac{1}{2}F''(x_n)e_n^2+O(\|e_n\|^3)\Big
),
\]
or
\begin{equation}\label{eq:6}
F'(x_n)^{-1}F(x_n)=e_n-\frac{1}{2}F'(x_n)^{-1}F''(x_n)e_n^2+O(\|e_n\|^3)
\end{equation}
From iterative process of (MTN) (\ref{eq:4}) we have:
\begin{equation}\label{eq:7}
\Big[F'(x_n)+2F'\Big(\frac{x_n+y_n}{2}\Big)+F'(y_n)\Big]e_{n+1}=\Big[F'(x_n)+2F'\Big(\frac{x_n+y_n}{2}\Big)+F'(y_n)\Big]e_n-4F(x_n).
\end{equation}
To continue we need the Taylor's expansion of $ F'(x_n-\theta
F'(x_n)^{-1}F(x_n))e_n$ as following:
\[
\begin{array}{ll}
F'(x_n-\theta F'(x_n)^{-1}F(x_n))e_n=F'(x_n)e_n-\theta
F''(x_n)F'(x_n)^{-1}F(x_n)e_n\hspace{1.5cm}\\
\hspace{5cm}+\frac{1}{2}\theta^2F^{(3)}(x_n)(F'(x_n)^{-1}F(x_n))^2e_n+O(\|e_n\|^4)
\end{array}
\]
by using (\ref{eq:6}) in above equation, we can write:
\[
\begin{array}{ll}
F'(x_n-\theta F'(x_n)^{-1}F(x_n))e_n\\
\hspace{2cm}=F'(x_n)e_n-\theta F''(x_n)\Big(e_n-\frac{1}{2}F'(x_n)^{-1}F''(x_n)e_n^2+O(\|e_n\|^3)\Big)e_n\\
\hspace{2.1cm}+\frac{1}{2}\theta^2F^{(3)}(x_n)\Big(e_n-\frac{1}{2}F'(x_n)^{-1}F''(x_n)e_n^2+O(\|e_n\|^3)\Big)^2e_n+O(\|e_n\|^4)
\end{array}
\]

after some manipulation we obtain:
\begin{equation}\label{eq:8}
\begin{array}{ll}
F'(x_n-\theta F'(x_n)^{-1}F(x_n))e_n\\
\hspace{2cm}=F(x_n)^{-1}e_n-\theta
F''(x_n)e_n^2+\frac{\theta}{2}F''(x_n)F'(x_n)^{-1}F''(x_n)e_n^3\\
\hspace{2.3cm}+\frac{\theta^2}{2}F^{(3)}(x_n)e_n^3+O(\|e_n\|^4),
\end{array}
\end{equation}
using (\ref{eq:8}) for $\theta=\frac{1}{2}$ and $\theta=1$, it is
obtained, respectively:
\begin{equation}\label{eq:9}
\begin{array}{ll}
 F'\Big(\frac{x_n+y_n}{2}\Big)=F'(x_n-\frac{1}{2}F'(x_n)^{-1}F(x_n))e_n\\
 \hspace{2.06cm}=F'(x_n)e_n-\frac{1}{2}F''(x_n)e_n^2+\frac{1}{4}F''(x_n)F'(x_n)^{-1}F''(x_n)e_n^3\\
 \hspace{2.2cm}+\frac{1}{8}F^{(3)}(x_n)e_n^3+O(\|e_n\|^4),
 \end{array}
 \end{equation}
 and(for $\theta=1$)
 \begin{equation}\label{eq:10}
 \begin{array}{ll}
F'(y_n)=F'(x_n-F'(x_n)F(x_n))e_n=F'(x_n)e_n-F''(x_n)e_n^2\\
\hspace{1.8cm}+\frac{1}{2}F''(x_n)F'(x_n)^{-1}F''(x_n)e_n^3+\frac{1}{2}F^{(3)}(x_n)e_n^3+O(\|e_n\|^4).
\end{array}
\end{equation}
by using (\ref{eq:5}),(\ref{eq:9}) and (\ref{eq:10}), we can write
right hand of Eq. (\ref{eq:7}) as following:
\[
\begin{array}{ll}
\Big(F'(x_n)+2F'(\frac{x_n+y_n}{2})+F'(y_n)\Big)e_n-4F(x_n)\\
\hspace{1cm}=\Bigg\{
F'(x_n)e_n+2\Big(F'(x_n)e_n-\frac{1}{2}F''(x_n)e_n^2\\
\hspace{1.5cm}+\frac{1}{4}F''(x_n)F'(x_n)^{-1}F''(x_n)e_n^3+\frac{1}{8}F^{(3)}(x_n)e_n^3+O(\|e_n\|^4)\Big)\\
\hspace{1.3cm}+\Big(F'(x_n)e_n-F''(x_n)e_n^2
+\frac{1}{2}F''(x_n)F'(x_n)^{-1}F''(x_n)e_n^3\\
\hspace{1.5cm}+\frac{1}{2}F^{(3)}(x_n)e_n^3+O(\|e_n\|^4)\Big)\Bigg\}-4\Bigg\{F'(x_n)e_n-\frac{1}{2}F''(x_n)e_n^2+\frac{1}{3!}F^{(3)}e_n^3\Bigg\}\\
\hspace{1cm}=F''(x_n)F'(x_n)^{-1}F''(x_n)e_n^3-\frac{1}{24}F^{(3)}(x_n)e_n^3+O(\|e_n\|^4)\\
\hspace{1cm}=\Big(
F''(x_n)F'(x_n)^{-1}F''(x_n)-\frac{1}{24}F^{(3)}(x_n)\Big)e_n^3+O(\|e_n\|^4)
\end{array}
\]
so using deduced  result, from Eq. (\ref{eq:7}) we obtain:
\[
\begin{array}{ll}
\Big(F'(x_n)+2F'(\frac{x_n+y_n}{2})+F'(y_n)\Big)e_{n+1}\\
\hspace{2cm}=\Big(
F''(x_n)F'(x_n)^{-1}F''(x_n)-\frac{1}{24}F^{(3)}(x_n)\Big)e_n^3+O(\|e_n\|^4).
\end{array}
\]
this  prove that the order of convergence is three, then the proof
is complete.$ \hspace{.5cm}\Box$
\section{Numerical examples}
In this section we will check the effectiveness of MTN (\ref{eq:4})
and other iterative methods in section \ref{sec:2}. All computations
are done by using $\bf{mathematica}$, stopping criteria
$\|x_{n+1}-x_n\|+\|F(x_n)\|\leqslant\epsilon$ is used for computer
programs. We use $\epsilon\leqslant 10^{-14}$.\\

\[
\begin{array}{ll}

 (a)\left\{\begin{array}{ll}
  e^{-x_1}. e^{x_2}+x_1 \cos(x_2)=0\\
  x_1+x_2-1=0.
\end{array}\right.
 &\qquad(b) \left\{\begin{array}{ll}
 x_1^2+3\ln(x_1)-x_2^2=0\\
 2x_1^2-x_1x_2-5x_1+1=0
 \end{array}\right.     \\
 \\
 (c) \left\{\begin{array}{ll}
x_1+2x_2-3=0\\
2x_1^2+x_2^2-5=0
\end{array}\right.
&\qquad(d)  \left\{\begin{array}{ll} \ln(x_1^2)-2\ln(\cos(x_2))=0\\
x_1\tan(\frac{x_1}{\sqrt{2}}+x_2)-\sqrt{2}=0
\end{array}\right.   \\
\\
  (e)\left\{\begin{array}{ll}
 x_1+e^{x_2}-\cos(x_2)=0\\
3x_1-x_2-\sin(x_2)=0
\end{array}\right.
&\qquad(f)  \left\{\begin{array}{ll}
x_1^2+x_2^2+x_3^2=9\\
 x_1.x_2.x_3-1=0\\
x_1+x_2-x_3^2=0
\end{array}\right.  \\
\\
  (g)\left\{\begin{array}{ll}
\cos(x_2)-\sin(x_1)=0\\
(x_3)^{x_1}-\frac{1}{x_2}=0\\
e^{x_1}-x_3^2=0
\end{array}\right.
&\qquad(h)  \left\{\begin{array}{ll}
x_2x_3+x_4(x_2+x_3)=0\\
 x_1x_3+x_4(x_1+x_3)=0\\
x_1x_2+x_4(x_1+x_2)=0 \\
x_1x_2+x_1x_3+x_2x_3=1
\end{array}\right.
\end{array}
\]
 \hspace{-1.cm}\small Approximations of $x_i$s for examples (a)-(e).\hspace{1.5cm}$\bf{Table\,\, 1}$

\[\hspace{-1.7cm}
 \begin{array}{|cccccc|}

\hline\
 F(x)& x_0& Method & Approximated\, solution & Iteration & Error\,
 estimation\\
 \hline\
 (a)&(1,2)& CN& (46.61144449,-45.61144449) & 9& 9.33\times
 10^{-15}\\
 &&TN&(-4.38161975,5.38161976)&5&9.10\times 10^{-14}\\
 &&MN&(-12.92527753,13.92527753)&6&3.55\times10^{-15}\\
 &&HN&&&\\
 &&TMN&-16.44481890,17.44481890&6&1.64\times 10^{-14}\\
 \hline\
 (b)&(3.4,2.2)&CN& (5.26375932,5.71748439)&10&7.99\times 10^{-15}\\
 &&TN&(5.26375932,5.71748439)&7&6.21\times 10^{-15}\\
  &&MN&(5.26375932,5.71748439)&7&4.44\times 10^{-15}\\
   &&HN&(5.26375932,5.71748439)&8&4.44\times 10^{-15}\\
   &&MTN&(5.26375932,5.71748439)&8&8.88\times 10^{-15}\\
   \hline\
   (c)&(1.5,1)&CN&(1.48803387,0.7558306)&5&9.99\times 10^{-16}\\
   &&TN&(1.48803387,0.7558306)&4&9.99\times 10^{-16}\\
   &&MN&(1.48803387,0.7558306)&4&9.99\times 10^{-16}\\
   &&HN&(1.48803387,0.7558306)&5&9.99\times 10^{-15}\\
   &&MTN&(1.48803387,0.7558306)&4&9.99\times 10^{-16}\\
   \hline\
   (d)&(0.2,0.2)&CN&No\, \,convergence&-&-\\
   &&TN&No\,\,convergence&-&-\\
   &&MN&No\,\,convergence&-&-\\
   &&HN&No\,\,convergence&-&-\\
   &&MTN&(0.95480414,6.58498148)&8&1.33\times 10^{-15}\\
   \hline\
   (e)&(-1,-3)&CN&(2.1378\times10^{-16},3.207\times10^{-16})&457&6.22\times10^{-16}\\
   &&TN&(4.817\times10^{-18},7.221\times10^{-18})&28&1.21\times10^{-27}\\
   &&MN&(-2.471\times10^{-17},-3.707\times10^{-27})&8&4.85\times10^{-2}\\
&&HN&(6.715\times10^{-17},9.107\times10^{-17})&53&2.77\times10^{-22}\\
&&MTN&(-1.364\times10^{-17},-2.046\times1)^{-17})&49&0\\
\hline
\end{array}
\]
 \hspace{-1.5cm}\small Approximations of $x_i$s for examples (f)-(h).\hspace{1.5cm}$\bf{Table\,\, 2}$

\[\hspace{-1.7cm}
 \begin{array}{|cccccc|}

\hline\
 F(x)& x_0& Method & Approximated\, solution & Iteration & Error\,
 estimation\\
 \hline\

 (f)&(2,2,0.5)&CN&(-2.090295,2.140258,-0.223525)&8&8.88\times10^{-16}\\
&&TN&(-2.090295,2.140258,-0.223525)&5&8.88\times10^{-16}\\
&&MN&(-2.090295,2.140258,-0.223525)&5&9.02\times10^{-16}\\
&&HN&(-2.090295,2.140258,-0.223525)&6&1.78\times10^{-15}\\
&&MTN&(-2.090295,2.140258,-0.223525)&5&9.02\times10^{-16}\\
\hline\
 (g)&(-2.5,1,1)&CN&(0.909569,0.661227,1.575834)&10&6.82\times10^{-14}\\
 &&TN&No\, convergence&-&-\\
 &&MN&(0.909569,0.661227,1.575834)&5&8.48\times10^{-14}\\
 &&HN&No\,\,convergence&-&-\\
 &&MTN&No\,\,convergence&-&-\\
 \hline\
(h)&(0.5,0.5,0.5,0.2)&CN&(0.5773,0.5773,0.5773,-0.2886)&5&2.22\times10^{-16}\\
&&TN&(0.5773,0.5773,0.5773,-0.2886)&4&1.11\times10^{-16}\\
&&MN&(0.5773,0.5773,0.5773,-0.2886)&4&1.11\times10^{-16}\\
&&HN&(0.5773,0.5773,0.5773,-0.2886)&6&1.31\times10^{-13}\\
&&MTN&(0.5773,0.5773,0.5773,-0.2886)&4&1.11\times10^{-16}\\
\hline
 \end{array}
 \]


\end{document}